\documentclass[12pt,letterpaper]{amsart}

\usepackage{xcolor}
\usepackage[english]{babel}
\usepackage[T1]{fontenc}
\usepackage[utf8]{inputenc}

\usepackage{amsmath,amssymb,amsthm,mathrsfs}
\usepackage{tikz-cd}
\usepackage{graphicx}
\usepackage{xypic}
\usepackage{hyperref}

\hypersetup{
  colorlinks=true,
  linkcolor=blue,
  urlcolor=red,
  pdftitle={}
}

\newtheorem{theorem}{Theorem}[section]
\newtheorem{proposition}[theorem]{Proposition}

\newtheorem{lemma}[theorem]{Lemma}

\newtheorem{definition}[theorem]{Definition}

\title{On Amoebas and the Arithmetic Rigidity of the Mahler Measure}
\author{Mounir Hajli}



\address{Westlake university}

\email{hajlimounir@gmail.com}

\date{\today}

\begin{document}

\maketitle

\begin{abstract}

We develop a geometric rigidity theory for Mahler measure based on the global structure of Ronkin functions and amoebas. 
Our starting point is the theorem of Ronkin  describing the affine linearity of the Ronkin function on the complement of the amoeba and the identification of its gradients with lattice points of the Newton polytope. 
Combining this deep analytic input with convex duality and toric  geometry, we show that the vanishing of the Mahler measure is an extremal phenomenon: it occurs if and only if the Ronkin function collapses to the piecewise linear support function of the Newton polytope.

This collapse forces a degeneration of the amoeba onto the $(n\!-\!1)$-skeleton of a toric fan and implies that the zero locus of the polynomial is contained in a translate of a codimension-one subtorus. 
As a first application of this rigidity principle, we recover the generalized Kronecker theorem characterizing Laurent polynomials of vanishing Mahler measure. 
Our approach bypasses the classical Boyd--Lawton limit formulas and replaces analytic projection arguments by a global convex-geometric and potential-theoretic mechanism rooted in the geometry of amoebas.

\end{abstract}

\tableofcontents

\section{Introduction}

The Mahler measure of a Laurent polynomial 
$f \in \mathbb{Z}[x_1^{\pm1},\dots,x_n^{\pm1}]$, defined as the integral of $\log|f|$ over the compact torus $\mathbb{S}^n$, is a fundamental invariant linking Diophantine geometry, the distribution of algebraic numbers on the unit torus, and the global geometry of hypersurfaces in algebraic tori. 
Since the pioneering work of Lehmer and Mahler \cite{Lehmer1933, Mahler}, a central problem has been to understand the arithmetic meaning of small or vanishing Mahler measure and to characterize those polynomials for which $m(f)=0$. 
In the univariate case this is answered by Kronecker's theorem, while in several variables it was resolved through the work of Boyd, Lawton, and Smyth using delicate analytic methods based on Fourier analysis and limiting procedures.

A different geometric structure underlying Mahler measure emerged from the theory of amoebas and Ronkin functions. 
To a Laurent polynomial $f$ one associates its amoeba $\mathcal{A}_f \subset \mathbb{R}^n$ and its Ronkin function $N_f$, a convex potential obtained by averaging $\log|f|$ over real tori. 
A fundamental theorem of Ronkin asserts that $N_f$ is affine linear on each connected component of the complement $\mathbb{R}^n \setminus \mathcal{A}_f$, and that the gradient on such a component is a lattice point of the Newton polytope of $f$. 
This result, which reflects deep analytic and pluripotential-theoretic properties of $f$, provides a precise link between the geometry of the amoeba and the combinatorics of the Newton polytope, and plays a decisive role in the present work.

\medskip

In this paper, we propose a different approach to the vanishing of Mahler measure, rooted in the global geometry of amoebas and Ronkin functions and organized using the language of toric Arakelov geometry. 
Rather than relying on Boyd--Lawton type limit formulas or analytic projection arguments, we interpret the Mahler measure as a height only after identifying a rigid geometric mechanism governing the associated convex potential. 
From this perspective, the condition $m(f)=0$ emerges as an extremal phenomenon reflecting a collapse of the Ronkin function to a purely combinatorial object.

Our proof relies on the interplay between three key ingredients:
\begin{enumerate}
    \item \emph{The Ronkin function.}  
    To a Laurent polynomial $f$ one associates a convex function $N_f$ on $N_{\mathbb{R}}$, defined as the average of $\log|f|$ over real tori.  
    A fundamental theorem of Ronkin  asserts that $N_f$ is affine linear on each connected component of the complement $N_{\mathbb{R}}\setminus\mathcal{A}_f$, and that the corresponding gradients are lattice points of the Newton polytope.  
    This deep analytic result, which reflects subtle pluripotential-theoretic properties of $f$, is a decisive input in our argument.
    
    \item \emph{The amoeba $\mathcal{A}_f$.}  
    The amoeba, defined as the image of the zero locus of $f$ under the logarithmic map, records precisely where the Ronkin function fails to be affine.  
    It acts as the singular support of the Monge--Amp\`ere measure associated with $N_f$ and controls the regions where new gradients appear.  
    The geometry of $\mathcal{A}_f$ therefore governs the possible deviations of $N_f$ from a piecewise linear support function.

    \item \emph{Convex duality and toric metrics.}  
    The Legendre--Fenchel transform of $N_f$ produces a concave function on the Newton polytope, which we interpret as a toric metric on the associated divisor.  
    In this step, convex duality plays a bookkeeping role, allowing one to globalize the Ronkin potential and relate its total Monge--Amp\`ere mass to an arithmetic height.
\end{enumerate}

We show that for $f \in \mathbb{Z}[M]$, the condition $m(f)=0$ forces the Ronkin function to coincide identically with the support function of its Newton polytope. 
This rigidity phenomenon means that no new gradients occur beyond those dictated combinatorially by $\mathrm{NP}(f)$, and hence that the Ronkin--amoeba structure theorem degenerates to its extremal case. 
Geometrically, this forces a collapse of the amoeba onto the $(n\!-\!1)$-skeleton of a toric fan, confining $\mathcal{A}_f$ to a union of affine hyperplanes determined by the faces of $\mathrm{NP}(f)$.  
This degeneration implies that the zero locus of $f$ is contained in a translate of a codimension-one subtorus, forcing $f$ to be a binomial.  
The cyclotomic property then follows from Galois equivariance.

This methodology bypasses limit formulas entirely and treats the Mahler measure as a global geometric invariant extracted from the extremal behavior of the Ronkin potential. 
Beyond yielding a new proof of the generalized Kronecker theorem, this framework suggests that the Lehmer problem, the search for a positive lower bound for nonzero Mahler measures, may be reformulated as a stability problem for Monge--Amp\`ere measures associated with Ronkin functions on toric varieties.

\section{Mahler Measures and the Kronecker Theorem}

Let $M \cong \mathbb{Z}^n$ be a lattice of rank $n$, and let $\mathbb{C}[M]$ denote the algebra of Laurent polynomials with complex coefficients. For an element $f = \sum_{m \in M} c_m \chi^m \in \mathbb{C}[M]$, where $\chi^m$ denotes the character associated with the lattice point $m$, the (logarithmic) Mahler measure is defined as the integral of $\log |f|$ over the compact torus $\mathbb{S}^n \subset (\mathbb{C}^\times)^n$ with respect to the normalized Haar measure $d\mu$:
\begin{equation}
m(f) = \int_{\mathbb{S}^n} \log |f(z)| \, d\mu(z) = \frac{1}{(2\pi i)^n} \int_{|z_1|=1} \cdots \int_{|z_n|=1} \log |f(z_1, \ldots, z_n)| \frac{dz_1 \cdots dz_n}{z_1 \cdots z_n}.
\end{equation}

In the univariate case, for a polynomial $f \in \mathbb{C}[x]$ with $f(x) = a \prod_{i=1}^d (x - \alpha_i)$, Jensen's formula provides the explicit evaluation:
\begin{equation}\label{mahler1}
m(f) = \log |a| + \sum_{i=1}^d \log^+ |\alpha_i|,
\end{equation}
where $\log^+ t = \log \max(1, t)$. For a polynomial $f \in \mathbb{Z}[x]$, it follows immediately from \eqref{mahler1} that $m(f) \geq 0$. A classical theorem of Kronecker characterizes the vanishing of this measure: for $f \in \mathbb{Z}[x]$, $m(f) = 0$ if and only if $f$ is a product of a monomial and cyclotomic polynomials.

This result was generalized to the multivariate setting by Boyd \cite{Boyd}, Lawton \cite{Lawton-Kronecker}, and Smyth \cite{Smyth}. We recall that a Laurent polynomial $f \in \mathbb{Z}[M]$ is defined as a \emph{generalized cyclotomic polynomial} if it is of the form
\begin{equation}
f = \chi^\lambda \Phi_k(\chi^\mu)
\end{equation}
for some multi-indices $\lambda, \mu \in \mathbb{Z}^n$ and a cyclotomic polynomial $\Phi_k \in \mathbb{Z}[t]$.

\begin{theorem}[Boyd-Lawton-Smyth]\label{Kronecker}
Let $f \in \mathbb{Z}[M]$ be an irreducible Laurent polynomial. Then $m(f) = 0$ if and only if $f$ is a generalized cyclotomic polynomial.
\end{theorem}

The objective of the present work is to demonstrate that this vanishing condition is a manifestation of a deeper arithmetic rigidity within the framework of Arakelov geometry.

\section{Convex Geometry and Legendre--Fenchel Duality}

In this section, we establish the necessary framework from convex analysis, following the conventions suitable for toric Arakelov geometry as developed in \cite{Burgos2, Oda, Roc70}.

Let $N \cong \mathbb{Z}^n$ be a free abelian group of rank $n$, and let $M = \mathrm{Hom}_{\mathbb{Z}}(N, \mathbb{Z})$ be its dual lattice. We denote by $N_{\mathbb{R}} = N \otimes_{\mathbb{Z}} \mathbb{R}$ and $M_{\mathbb{R}} = M \otimes_{\mathbb{Z}} \mathbb{R}$ the corresponding real vector spaces, equipped with the canonical pairing $\langle \cdot, \cdot \rangle: M_{\mathbb{R}} \times N_{\mathbb{R}} \to \mathbb{R}$.

\begin{theorem}[\cite{Oda}, Theorem A.18]\label{Minkowski_Duality}
Let $\mathcal{K}(M_{\mathbb{R}})$ denote the set of non-empty compact convex subsets of $M_{\mathbb{R}}$, and let $\mathrm{SF}(N_{\mathbb{R}})$ denote the set of positively homogeneous concave functions on $N_{\mathbb{R}}$. There exists a bijective correspondence between $\mathcal{K}(M_{\mathbb{R}})$ and $\mathrm{SF}(N_{\mathbb{R}})$ given by the mutually inverse mappings $K \mapsto \psi_K$ and $\psi \mapsto K_\psi$, defined by:
\begin{enumerate}
    \item $\psi_K(u) = \inf_{v \in K} \langle v, u \rangle$ for $u \in N_{\mathbb{R}}$;
    \item $K_\psi = \{ v \in M_{\mathbb{R}} \mid \langle v, u \rangle \geq \psi(u), \forall u \in N_{\mathbb{R}} \}$.
\end{enumerate}
Furthermore, this correspondence is additive with respect to the Minkowski sum in $\mathcal{K}(M_{\mathbb{R}})$ and the pointwise sum in $\mathrm{SF}(N_{\mathbb{R}})$.
\end{theorem}

Crucial to our development is the Legendre--Fenchel transform, which allows us to pass from potentials on $N_{\mathbb{R}}$ to weights on the Newton polytope. 

\begin{definition}\label{def:Legendre-Fenchel}
Let $f: N_{\mathbb{R}} \to \mathbb{R}$ be a concave function. The \emph{stability set} of $f$, denoted $\mathrm{stab}(f)$, is the set of $x \in M_{\mathbb{R}}$ such that $u \mapsto \langle x, u \rangle - f(u)$ is bounded from below. For a concave function with $\mathrm{stab}(f) \neq \emptyset$, the \emph{Legendre--Fenchel dual} $f^\vee: \mathrm{stab}(f) \to \mathbb{R}$ is defined by:
\[
f^\vee(x) = \inf_{u \in N_{\mathbb{R}}} \left( \langle x, u \rangle - f(u) \right).
\]
\end{definition}

In the context of Mahler measures, the stability set of the Ronkin function $\rho_f$ will be precisely the Newton polytope $\mathrm{NP}(f)$, and $f^\vee$ will represent the \emph{roof function} describing the canonical metric on the associated toric line bundle.

\section{Toric Varieties and Equivariant Divisors}

We briefly summarize the dictionary between the combinatorial data of fans and the geometry of toric varieties. For comprehensive treatments, we refer the reader to \cite{Burgos2, Fulton, Oda}.

\subsection{Cones, Fans, and Varieties}
Let $N \cong \mathbb{Z}^n$ be a lattice and $M$ its dual. A \emph{fan} $\Sigma$ in $N_{\mathbb{R}}$ is a finite collection of strongly convex rational polyhedral cones such that every face of a cone in $\Sigma$ is also in $\Sigma$, and the intersection of any two cones is a face of each. We denote by $X_\Sigma$ the associated $n$-dimensional toric variety, and $T_N$ its associated torus acting on $X_\Sigma$. We assume throughout that $\Sigma$ is complete, i.e., its support $|\Sigma| = \bigcup_{\sigma \in \Sigma} \sigma$ is the entire space $N_{\mathbb{R}}$, which implies that $X_\Sigma$ is proper over $\mathbb{C}$.

\subsection{Equivariant Divisors and Support Functions}
Let $\mathrm{SF}(N, \Sigma)$ denote the group of $\Sigma$-linear support functions $h: N_{\mathbb{R}} \to \mathbb{R}$. These functions are uniquely determined by a set of lattice points $\{l_\sigma\}_{\sigma \in \Sigma_{\max}} \subset M$ such that $h(n) = \langle l_\sigma, n \rangle$ for all $n \in \sigma$. Such a function $h$ corresponds to a $T_N$-invariant Cartier divisor $D_h$.

The properties of the line bundle $\mathcal{O}(D_h)$ are encoded in the convexity of $h$. We define the \emph{associated polytope} $\square_h \subset M_{\mathbb{R}}$ as:
\[
\square_h := \{ m \in M_{\mathbb{R}} \mid \langle m, n \rangle \geq h(n), \, \forall n \in N_{\mathbb{R}} \}.
\]

The following characterization is central to our height-theoretic approach:

\begin{proposition}[\cite{Oda}, Prop 2.1, Thm 2.13]\label{prop:toric_correspondence}
Let $h \in \mathrm{SF}(N, \Sigma)$. The following are equivalent:
\begin{enumerate}
    \item The equivariant divisor $D_h$ is ample.
    \item $h$ is strictly concave with respect to $\Sigma$, i.e., $h(n+n') \geq h(n) + h(n')$ with equality if and only if $n, n'$ lie in the same cone $\sigma \in \Sigma$.
    \item The polytope $\square_h$ is $n$-dimensional and the mapping $\sigma \mapsto \sigma^\dagger$ gives a bijection between the cones of $\Sigma$ and the faces of $\square_h$, with $\dim \sigma + \dim \sigma^\dagger = n$.
\end{enumerate}
\end{proposition}

Under these conditions, the vertices of $\square_h$ are precisely the set $\{l_\sigma \mid \sigma \in \Sigma_{\max}\}$.

\subsection{Toric metrized divisors and roof functions}\label{subsec:toric-metrics}

In this subsection, we recall some standard facts from toric Arakelov geometry (\cite{Burgos3,Burgos2}) that will be used as a formal framework to organize the convex-analytic objects arising later from Ronkin functions.  
Our main rigidity argument does \emph{not} originate in this formalism: the decisive analytic input will come from structural results on Ronkin functions and amoebas.  
The role of toric metrized divisors is to provide a global geometric language in which these convex potentials can be interpreted and their total Monge--Amp\`ere mass identified with an arithmetic height.

Let $X = X_\Sigma$ be a proper toric variety over $\mathbb{C}$ with dense torus $T_N$.

\begin{definition}\label{def:toric-metric}
A \emph{metrized divisor} $\overline{D} = (D, \|\cdot\|)$ consists of a Cartier divisor $D$ together with a continuous metric on the line bundle $\mathcal{O}(D)$.  
The metric is called \emph{toric} if $D$ is $T_N$-invariant and the metric is invariant under the action of the compact torus $\mathbb{S}^n \subset T_N$.  

To such a metric one associates a function $\psi_{\overline{D}}: N_{\mathbb{R}} \to \mathbb{R}$ defined by
\[
\psi_{\overline{D}}(u) = \log \|s_D(z)\|,
\]
where $z \in T_N(\mathbb{C})$ is any point such that $-\mathrm{Log}(z) = u$, and $s_D$ denotes the canonical toric section of $\mathcal{O}(D)$.
\end{definition}

The function $\psi_{\overline{D}}$ is well defined and relates to the Green function $g_{\overline{D}}$ of the metrized divisor via
\[
\psi_{\overline{D}}(u) = - g_{\overline{D}}(z).
\]
A toric metric is \emph{semipositive} if and only if the associated function $\psi_{\overline{D}}$ is concave on $N_{\mathbb{R}}$.

\begin{definition}\label{def:roof-function}
Let $\overline{D}$ be a semipositive toric metrized divisor.  
The associated \emph{roof function} $\vartheta_{\overline{D}}: \square_D \to \mathbb{R}$ is the concave function on the Newton polytope $\square_D$ defined as the Legendre--Fenchel dual of $\psi_{\overline{D}}$:
\[
\vartheta_{\overline{D}}(x)
= \psi_{\overline{D}}^\vee(x)
= \inf_{u \in N_{\mathbb{R}}} \bigl( \langle x, u \rangle - \psi_{\overline{D}}(u) \bigr).
\]
\end{definition}

The following result summarizes the dictionary between toric metrics and convex geometry.

\begin{proposition}[\cite{Burgos3}, Prop.~4.19]\label{prop:toric_metrics}
Let $D$ be an ample toric divisor on $X$ with associated support function $\Psi_D$.  
There is a one-to-one correspondence between:
\begin{enumerate}
    \item Toric semipositive metrics $\|\cdot\|$ on $\mathcal{O}(D)$;
    \item Concave functions $\psi: N_{\mathbb{R}} \to \mathbb{R}$ such that $|\psi - \Psi_D|$ is bounded on $N_{\mathbb{R}}$;
    \item Continuous concave functions $\vartheta: \square_D \to \mathbb{R}$ on the Newton polytope $\square_D$.
\end{enumerate}
\end{proposition}

Under this correspondence, the \emph{canonical metric} $\|\cdot\|_{\mathrm{can}}$ (which has zero curvature on the dense torus) corresponds to the support function $\psi = \Psi_D$ and to the trivial roof function $\vartheta \equiv 0$.

\medskip

\noindent
\textbf{Remark.}
In the sequel, Ronkin functions of Laurent polynomials will be shown to define particular concave functions $\psi$ on $N_{\mathbb{R}}$ satisfying the boundedness condition of Proposition~\ref{prop:toric_metrics}.  
The purpose of this subsection is therefore not to introduce new tools, but to provide a geometric container in which the rigidity properties of Ronkin functions, coming from amoeba geometry and pluripotential theory, can be interpreted in terms of toric heights.

\section{The Ronkin function and the amoeba}\label{sec:ronkin}

The Ronkin function is the central analytic object in our approach.  
It provides a canonical extension of the Mahler measure to a global potential on $N_{\mathbb{R}}$ and encodes, in a highly rigid way, the interaction between the zero locus of a Laurent polynomial and the geometry of its Newton polytope.  
The decisive input in our proof of the Kronecker theorem will come from deep structural properties of this function established by Ronkin.

\subsection{Definitions and convexity}

Let $f \in \mathbb{C}[M]$ be a Laurent polynomial.  
The \emph{Ronkin function} $N_f : N_{\mathbb{R}} \to \mathbb{R}$ is defined by
\begin{equation}
N_f(u)
= \frac{1}{(2\pi i)^n}
\int_{\mathrm{Log}^{-1}(u)}
\log |f(z)|
\frac{dz_1 \wedge \cdots \wedge dz_n}{z_1 \cdots z_n}.
\end{equation}
Equivalently, $N_f(u)$ is the Mahler measure of the twisted polynomial
\[
f_u(z) = f(e^{u_1} z_1, \dots, e^{u_n} z_n).
\]
A fundamental result of Ronkin \cite{Ronkin2001} asserts that $N_f$ is a convex function on $N_{\mathbb{R}}$.

The \emph{amoeba} of $f$ is the subset
\[
\mathcal{A}_f = \mathrm{Log}(\mathcal{Z}(f)) \subset N_{\mathbb{R}},
\]
where $\mathcal{Z}(f)$ denotes the zero locus of $f$ in $(\mathbb{C}^\times)^n$.  
The large-scale geometry of $N_f$ is governed by the topology of $\mathcal{A}_f$, and in particular by the structure of the complement $N_{\mathbb{R}} \setminus \mathcal{A}_f$.

\begin{theorem}[\cite{PassareRullgard2004, Ronkin2001}]\label{Ronkin_Properties}
The Ronkin function $N_f$ is affine linear on each connected component of the complement
$N_{\mathbb{R}} \setminus \mathcal{A}_f$.  
More precisely, if $U$ is such a component, then the gradient $\nabla N_f|_U$ is constant and equal to a lattice point
\[
m \in M \cap \mathrm{NP}(f),
\]
called the \emph{order} of the component.
\end{theorem}

\medskip

\noindent
This theorem is a deep and highly nontrivial result, relying on complex analytic and pluripotential-theoretic techniques.  
It expresses a strong rigidity phenomenon: away from the amoeba, the Ronkin function is not merely smooth or harmonic, but \emph{exactly affine linear}, with slopes quantized by the lattice points of the Newton polytope.  
This piecewise-linear structure will be the key mechanism behind our rigidity argument for Mahler measure.

\subsection{The concave potential $\rho_f$}

For later geometric interpretation, it is convenient to introduce the concave function
\[
\rho_f(u) = - N_f(u).
\]
Let $\Psi_f$ denote the concave support function of the Newton polytope $\mathrm{NP}(f)$.

The following result, due to Passare--Rullg{\aa}rd and refined in later work, shows that $\rho_f$ is a global regularization of $\Psi_f$.

\begin{proposition}[\cite{PassareRullgard2004, Gualdi}]\label{bound}
The potential $\rho_f$ approximates the support function $\Psi_f$ in the following precise sense:
\begin{enumerate}
\item The difference $|\rho_f - \Psi_f|$ is bounded on $N_{\mathbb{R}}$.
\item For all $u \in N_{\mathbb{R}}$, one has the inequality
\begin{equation}\label{Passare_Ineq}
\rho_f(u)
\leq \Psi_f(u)
- \log \min_{m \in \mathcal{V}(f)} |c_m|,
\end{equation}
where $\mathcal{V}(f)$ denotes the set of vertices of $\mathrm{NP}(f)$.
\end{enumerate}
\end{proposition}

\medskip

\noindent
In particular, $\rho_f$ has the same asymptotic slopes as $\Psi_f$ at infinity, and its deviation from the support function is globally controlled by the extreme coefficients of $f$.  
This boundedness property will later allow us to interpret $\rho_f$ as a toric metric, but it plays no role in the rigidity phenomenon itself, which is entirely encoded in Theorem~\ref{Ronkin_Properties}.

\subsection{Toric realization of the Ronkin potential}\label{subsection:toric-realization}

To embed the Ronkin potential into the formalism of toric Arakelov geometry, one must ensure that the support function $\Psi_f$ is compatible with the fan of a projective toric variety.  
This requires a standard toric resolution procedure.

\begin{lemma}\label{lemma:toric_resolution}
There exists a smooth projective toric variety $X_{\Sigma}$ and a proper toric morphism
\[
\pi : X_{\Sigma} \longrightarrow \mathbb{P}^n
\]
such that:
\begin{enumerate}
\item $\pi$ restricts to the identity on the dense open torus;
\item The support function $\Psi_f$ is $\Sigma$-linear.
\end{enumerate}
\end{lemma}

\begin{proof}
See \cite{Gualdi}.  
By the toric Chow lemma and the theory of fan refinements (cf.\ \cite{Oda, Fulton}), one can refine the fan of $\mathbb{P}^n$ to a fan $\Sigma$ for which $\Psi_f$ becomes piecewise linear.  
After further refinement to ensure regularity and projectivity, the induced morphism $\pi$ has the required properties.
\end{proof}

By pulling back to $X_{\Sigma}$, Proposition~\ref{bound} implies that $|\rho_f - \Psi_f|$ is bounded.  
Consequently, $\rho_f$ defines a continuous semipositive toric metric on the line bundle $\mathcal{O}(D_{\Psi_f})$ over $X_{\Sigma}$.

\medskip

\noindent
\textbf{Remark.}
The purpose of this toric realization is purely interpretative: it allows one to view $\rho_f$ as a metrized divisor and to identify its total Monge--Amp\`ere mass with an arithmetic height.  
None of the rigidity results proved later depend on this formalism; they rely instead on the affine-linearity and quantization properties of the Ronkin function stated in Theorem~\ref{Ronkin_Properties}.

\section{A geometric proof of the generalized Kronecker theorem}\label{sec:kronecker}

In this section, we prove Theorem~\ref{Kronecker}.  
The conceptual strategy is the following:

\begin{itemize}
\item The vanishing of the Mahler measure forces the Ronkin function to become extremal among all convex potentials with the same asymptotic slopes.
\item By the rigidity theorem of Ronkin, this extremality propagates to a collapse of the amoeba.
\item This collapse constrains the zero locus of $f$ to lie in a binomial hypersurface, from which the cyclotomic form follows by Galois invariance.
\end{itemize}

The decisive step is therefore analytic and geometric in nature; toric Arakelov geometry intervenes only to convert extremality of the Ronkin potential into the scalar condition $m(f)=0$.

\medskip

Let $f \in \mathbb{Z}[M]$ be an irreducible Laurent polynomial.  
Up to multiplying by a monomial, we may assume that $0 \in \mathrm{NP}(f)$.

By a theorem of Maillot~\cite{Maillot}, the Mahler measure of $f$ coincides with the canonical height of the associated hypersurface:
\[
h_{\overline{\mathcal{O}(1)}_{\mathrm{can}}}(\mathrm{div}(f)) = m(f).
\]
This identity allows one to interpret $m(f)$ as the total Monge--Amp\`ere mass of the concave potential $\rho_f = -N_f$.

Let $\check{\rho}_f$ denote the Legendre--Fenchel dual of $\rho_f$.  
From Proposition~\ref{bound} and the integrality of the coefficients $c_m$, we obtain the global bounds
\begin{equation}\label{ineq:duality_bound}
0
\leq
\log \bigl( \min_{m \in \mathcal{V}(f)} |c_m| \bigr)
\leq
\check{\rho}_f(x)
\leq
-\rho_f(0)
=
m(f),
\qquad x \in \mathrm{NP}(f).
\end{equation}

\medskip

\noindent
Thus, the Mahler measure controls the vertical deviation of $\check{\rho}_f$ above the Newton polytope.

\begin{proposition}[Convex rigidity]\label{prop:rigidity}
If $m(f)=0$, then $\rho_f(u)=\Psi_f(u)$ for all $u \in N_{\mathbb{R}}$.
\end{proposition}

\begin{proof}
Assume that $m(f)=0$.  
By the global bounds \eqref{ineq:duality_bound}, this forces
\[
\check{\rho}_f(x) = 0
\qquad \text{for all } x \in \mathrm{NP}(f).
\]

On the other hand, the Legendre--Fenchel dual of the support function $\Psi_f$ is also identically zero on $\mathrm{NP}(f)$.

\medskip

\noindent
To legitimately compare $\rho_f$ and $\Psi_f$ via convex duality, we must ensure that both functions belong to the same duality class of concave functions with fixed asymptotic behavior.  
This is guaranteed by Proposition~\ref{bound}, which shows that $|\rho_f - \Psi_f|$ is bounded on $N_{\mathbb{R}}$, and by Lemma~\ref{lemma:toric_resolution}, which provides a refinement of the fan such that $\Psi_f$ becomes a $\Sigma$-linear support function.

\medskip

\noindent
After pulling back to the corresponding projective toric variety $X_\Sigma$, both $\rho_f$ and $\Psi_f$ define concave functions whose difference is bounded and whose Legendre--Fenchel duals are finite precisely on the same polytope $\mathrm{NP}(f)$.  
In this setting, the Legendre--Fenchel transform is an involution on concave functions with these asymptotics.

\medskip

\noindent
Since $\check{\rho}_f = \check{\Psi}_f = 0$ on $\mathrm{NP}(f)$, the involution property of the Legendre--Fenchel transform yields
\[
\rho_f = \Psi_f
\qquad \text{on } N_{\mathbb{R}}.
\]

\medskip

\noindent
This argument uses convex duality only as a formal mechanism to convert a scalar vanishing condition into a global extremality statement.  
The geometric consequences of this extremality are provided by the rigidity properties of the Ronkin function recalled in Theorem~\ref{Ronkin_Properties}.
\end{proof}

\medskip

We now examine the implications of Proposition~\ref{prop:rigidity} for the zero locus of $f$.  
Let $\mathcal{A}_f \subset N_{\mathbb{R}}$ denote the amoeba of $f$.

\medskip

\noindent
By Proposition~\ref{prop:rigidity}, we have $\rho_f = \Psi_f$ on $N_{\mathbb{R}}$.  
Recall that $\Psi_f$ is a piecewise linear concave function whose domains of linearity are precisely the maximal cones $\sigma \in \Sigma_{\max}$ of the fan associated with the Newton polytope.

\medskip

\noindent
On the other hand, Theorem~\ref{Ronkin_Properties} asserts that the Ronkin function $N_f = -\rho_f$ is affine linear on each connected component of the complement $N_{\mathbb{R}} \setminus \mathcal{A}_f$, with gradient equal to a lattice point $m \in M \cap \mathrm{NP}(f)$ (the order of the component).

\medskip

\noindent
Since $\rho_f = \Psi_f$ and $\Psi_f$ is affine linear on each cone $\sigma \in \Sigma_{\max}$ but not across different cones, it follows that $\rho_f$ can only be affine linear on open subsets contained entirely in a single cone of $\Sigma$.

\medskip

\noindent
Therefore, every connected component of $N_{\mathbb{R}} \setminus \mathcal{A}_f$ must be contained in the relative interior of some maximal cone $\sigma \in \Sigma_{\max}$.  
Equivalently, the amoeba must be contained in the locus where $\Psi_f$ fails to be differentiable, namely the union of the codimension--one faces separating adjacent cones:

\[
\mathcal{A}_f 
\;\subseteq\; 
\bigcup_{\substack{\sigma, \sigma' \in \Sigma_{\max} \\ \sigma \neq \sigma'}} 
(\sigma \cap \sigma').
\]

\medskip

\noindent
Since $\Sigma$ is a complete fan, this union is a polyhedral complex of dimension at most $n-1$.  
In particular, $\mathcal{A}_f$ is contained in a finite union of affine hyperplanes in $N_{\mathbb{R}}$.

\medskip

We now show that the topological constraint on the amoeba forces $f$ to be binomial.

\medskip

\noindent
Since $\mathcal{A}_f$ is contained in a finite union of affine hyperplanes, there exists a nonzero lattice vector $m \in M$ and a real constant $c$ such that
\[
\langle m , u \rangle = c
\qquad \text{for all } u \in \mathcal{A}_f.
\]

\medskip

\noindent
Equivalently, for every $z \in Z(f) \subset (\mathbb{C}^\times)^n$, we have
\[
\langle m , \mathrm{Log}(z) \rangle = c,
\]
which implies
\[
|z^m| = e^{c}.
\]

\medskip

\noindent
Thus the holomorphic function $z \mapsto z^m$ has constant modulus on the analytic hypersurface $Z(f)$.  
By the maximum modulus principle applied to the regular locus of $Z(f)$, it follows that $z^m$ must be constant on each irreducible component of $Z(f)$.  
Since $f$ is irreducible, we conclude that
\[
z^m \equiv u
\qquad \text{on } Z(f)
\]
for some constant $u \in \mathbb{C}$ with $|u| = e^{c}$.

\medskip

\noindent
Consequently, the irreducible polynomial $f$ divides the binomial $z^m - u$ in $\mathbb{C}[M]$.  
Since $f \in \mathbb{Z}[M]$, its zero locus is invariant under the action of $\mathrm{Gal}(\overline{\mathbb{Q}}/\mathbb{Q})$, and therefore all Galois conjugates of $u$ also have the same modulus.

\medskip

\noindent
The condition $m(f)=0$ implies that $|u|=1$.  
By Kronecker's univariate theorem, this forces $u$ to be a root of unity.

\medskip

\noindent
Thus $f$ divides $z^m - \zeta$ for some $\zeta \in \mu_\infty$.  
Since $f$ is irreducible in $\mathbb{Z}[M]$, it follows that
\[
f(z) = a z^\lambda \Phi_k(z^\mu),
\]
where $\Phi_k$ is a cyclotomic polynomial and $a = \pm 1$.  
This completes the proof.

\bibliographystyle{plain}

\bibliography{biblio}

\end{document}